  %
  %
  %
  %
  %


  \newcount\fontset
  \fontset=1
  \def\dualfont#1#2#3{\font#1=\ifnum\fontset=1 #2\else#3\fi}
  \dualfont\eightrm {cmr8} {cmr7}
  \dualfont\eightsl {cmsl8} {cmr7}
  \dualfont\eightit {cmti8} {cmti10}
  \dualfont\eightmi {cmmi8} {cmmi7}
  \dualfont\eightbf {cmbx8} {cmbx10}
  \dualfont\fivemi {cmmi5} {cmmi7}
  \dualfont\tensc {cmcsc10} {cmcsc10}
  \dualfont\eighttt {cmtt8} {cmtt10}
  \dualfont\tenss {cmss10} {cmr10}
  \dualfont\titlefont {cmbx12} {cmbx10}
  \dualfont\eightsymbol {cmsy8} {cmsy10}


  \def\vg#1{\ifx#1\null\null\else
    \ifx#1\ { }\else
    \ifx#1,,\else
    \ifx#1..\else
    \ifx#1;;\else
    \ifx#1::\else
    \ifx#1''\else
    \ifx#1--\else
    \ifx#1))\else
    { }#1\fi \fi \fi \fi \fi \fi \fi \fi \fi}
  \newcount\secno \secno=0
  \newcount\stno
  \def\goodbreak{\vskip0pt plus.1\vsize\penalty-250
    \vskip0pt plus-.1\vsize\bigskip}
  \outer\def\section#1{\stno=0
    \global\advance\secno by 1
    \goodbreak\vskip\parskip
    \message{\number\secno.#1}
    \noindent{\bf\number\secno.\enspace#1.\enspace}}
  \def\seqnumbering{\global\advance\stno by 1
    \number\secno.\number\stno}
  \def\label#1{\global\edef#1{\number\secno.\number\stno}}
  \def\sysstate#1#2#3{\medbreak\noindent {\bf\seqnumbering.\enspace
    #1.\enspace}{#2 #3\vskip 0pt}\medbreak}
  \def\state#1 #2\par{\sysstate{#1}{\sl}{#2}}
  \def\definition#1\par{\sysstate{Definition}{\rm}{#1}}
  \def\proof{\medbreak\noindent{\it Proof.\enspace}}
  \def\proofend{\ifmmode\eqno\square\else\hfill\square\medbreak\fi}
  \newcount\zitemno \zitemno=0
  \def\zitem{\global\advance\zitemno by 1 \smallskip
    \item{\ifcase\zitemno\or i\or ii\or iii\or iv\or v\or vi\or vii\or
    viii\or ix\or x\or xi\or xii\fi)}}
  \def\$#1{#1 $$$$ #1}


  \def\({\left(}
  \def\){\right)}
  \def\[{\left\Vert}
  \def\]{\right\Vert}
  \def\*{\otimes}
  \def\+{\oplus}
  \def\:{\colon}
  \def\<{\langle}
  \def\>{\rangle}
  \def\and{\hbox{\quad and \quad}}
  \def\arw{\rightarrow}

  \def\crossproduct{\hbox to 1.8ex{$\times \kern-.45ex\vrule height1.1ex
depth0pt width0.45truept$\hfill}}
  \def\cstar{$C^*$}
  \def\for#1{,\quad \hbox{ for }#1} 
  \def\inv{^{-1}}

  \def\square{\hbox{$\sqcap\!\!\!\!\sqcup$}}
  \def\stress#1{{\it #1}\/}

  \def\|{\Vert}
  \def\inv{^{-1}}


  \newcount\bibno \bibno=0
  \def\newbib#1{\global\advance\bibno by 1 \edef#1{\number\bibno}}
  \def\cite#1{{\rm[\bf #1\rm]}}
  \def\scite#1#2{[{\bf #1},{\it #2\/}]}
  \def\lcite#1{#1}
  \def\se#1#2#3#4{\def\a{#1}\def\b{#2}\ifx\a\b#3\else#4\fi}
  \def\setem#1{\se{#1}{}{}{, #1}}
  \def\index#1{\smallskip \item{[#1]}}
  \def\tit#1{``#1''}

  \def\zarticle#1, auth = #2,
   title = #3,
   journal = #4,
   year = #5,
   volume = #6,
   pages = #7,
   NULL#8 xyzzy {\index{#1} #2, \tit{#3}, {\sl #4\/} {\bf #6} (#5), #7.}

  \def\ztechreport#1,
    auth = #2,
    title = #3,
    institution = #4,
    year = #5,
    type = #6,
    note = #7,
    NULL#8
    xyzzy
    {\index{#1} #2, \tit{#3}\setem{#6}\setem{#4}\setem{#5}\setem{#7}.}

  \def\zunpublished#1,
    auth = #2,
    title = #3,
    institution = #4,
    year = #5,
    type = #6,
    note = #7,
    NULL#8
    xyzzy
    {\index{#1} #2, \tit{#3}\setem{#6}\setem{#4}\setem{#5}\setem{#7}.}

  \def\zbook#1,
    auth = #2,
    title = #3,
    publisher = #4,
    year = #5,
    volume = #6,
    series = #7,
    NULL#8
    xyzzy
    {\index{#1} #2, \tit{#3}\setem{#7}\se{#6}{}{}{ vol. #6}, #4, #5.}

  \def\zmasterthesis#1,
    auth = #2,
    title = #3,
    school = #4,
    year = #5,
    type = #6,
    NULL#7
    xyzzy
    {\index{#1} #2, \tit{#3}, #6, #4, #5.}

  \def\zinproceedings#1,
  auth = #2,
    title = #3,
    booktitle = #4,
    year = #5,
    pages = #6,
    organization = #7,
    note = #8,
    NULL#9
    xyzzy
    {\index{#1} #2, \tit{#3}\se{#4}{}{}{, In
      {\sl #4}}\setem{#5}\setem{#6}\setem{#7}\setem{#8}.}

  \def\zphdthesis#1,
    auth = #2,
    title = #3,
    school = #4,
    year = #5,
    type = #6,
    NULL#7
    xyzzy
    {\index{#1} #2, \tit{#3}, #6, #4, #5.}

  \def\zbooklet#1,
    auth = #2,
    title = #3,
    howpublished = #4,
    year = #5,
    NULL#6
    xyzzy
    {\index{#1} #2, \tit{#3}, #4, #5.}

  \def\zmisc#1,
    auth = #2,
    title = #3,
    note = #4,
    howpublished = #5,
    year = #6,
    NULL#7
    xyzzy
    {\index{#1} #2, \tit{#3}\setem{#4}\setem{#5}\setem{#6}.}


  \nopagenumbers
  \voffset=2\baselineskip
  \advance\vsize by -\voffset
  \headline{\ifnum\pageno=1 \hfil \else \tensc\hfil
   partial representations and fell bundles
  \hfil\folio \fi}


  \def\assoc{\sim}
  \def\sassoc{\buildrel s \over \assoc}
  \def\pr{\sigma}
  \def\M{{\cal M}}
  \def\Lin{{\cal B}}
  \def\K{{\cal K}}
  \def\C{{\bf C}}
  \def\R{{\bf R}}
  \def\N{{\bf N}}
  \def\I{{\cal I}}
  \def\e{{\rm e}}
  \def\prep{partial representation\vg}
  \def\preps{partial representations\vg}
  \def\ap{approximation property\vg}
  \def\ss{semi-saturated\vg}
  \def\F{I\!\!F} 
  \def\B{I\!\!B}
  \def\Gen{{\cal S}}
  \def\a{\alpha}
  \def\b{\beta}
  \def\Pos{{\cal W}}
  \def\root{^{1\over 2}}
  \def\limn{\lim_{n\rightarrow\infty}}
  \def\=#1{\buildrel #1 \over =}
  \def\coef#1{\({n - |#1| + 1 \over n}\)\root}
  \def\ninv{\({1 \over n}\)\root}
  \def\pef#1{\Big( \e(#1) - f(#1) \Big)}
  \def\Cs{C^*(\B)}
  \def\Csr{C^*_r(\B)}


  \def\Claire{A}
  \def\CK{CK}
  \def\BGR{BGR}
  \def\TPA{E1}
  \def\Inverse{E2}
  \def\Amena{E3}
  \def\FD{FD}
  \def\Greenleaf{G}
  \def\Ng{Ng}
  \def\Nica{Ni}
  \def\Ped{P}
  \def\Nandor{S}
  \def\Zettl{Z}


  \null
  \centerline{\titlefont PARTIAL REPRESENTATIONS AND AMENABLE}
  \medskip
  \centerline{\titlefont FELL BUNDLES OVER FREE GROUPS}
  \bigskip
  \centerline{\tensc Ruy Exel\footnote{*}{\eightrm Partially supported
                     by CNPq -- Brazil.}}
  \bigskip
  \centerline{\bf May 25, 1997}

  \def\.#1{\hbox{\eightmi #1}}
  \def\II#1{\.I\kern-2pt\.{#1}}

  \bigskip
  \midinsert\narrower\narrower
  \baselineskip=11pt\eightbf
  \noindent ABSTRACT. We show that a Fell bundle
  $\II B =\{\.B_t\}_{t\in \F}$,
  over an arbitrary free group $\II F$, is amenable, whenever it is
orthogonal (in the sense that $\.B_x^*\.B_y=0$, if $\.x$ and $\.y$ are
distinct generators of $\II F$) and \ss (in the sense that $\.B_{ts}$
coincides with the closed linear span of $\.B_t\.B_s$, when the
multiplication ``$\.{ts}$'' involves no cancelation).
  \endinsert


  In this work we continue the study of the phenomena of amenability for
Fell bundles over discrete groups, initiated in \cite{\Amena}.  By
definition, a Fell bundle is said to be amenable if the left regular
representation of its cross-sectional \cstar-algebra is faithful.  This
property is also equivalent to the faithfulness of the standard
conditional expectation. The reader is referred to \cite{\Amena} for
more information, but we also offer
  a very brief survey containing some of the most relevant definitions,
 in our section on preliminaries below.

  The starting point for our work is Theorem 6.7 of \cite{\Amena}, where
it is shown that a certain grading of the Cuntz--Krieger algebra gives
rise to an amenable Fell bundle over a free group.  Our main goal is to
further pursue the argument leading to this result, in order to obtain a
large class of amenable Fell bundles.  We find that the crucial
properties implying the amenability of a Fell bundle,
  over a free group $\F$, are \stress{orthogonality} and
\stress{semi-saturatedness}.  A Fell bundle $\B=\{B_t\}_{t\in\F}$ is
said to be \stress{orthogonal} if the fibers $B_x$ and $B_y$,
corresponding to two distinct generators $x$ and $y$ of $\F$, are
orthogonal in the sense that $B_x^*B_y=0$.  On the other hand, $\B$ is
said to be \stress{\ss{}} when each fiber $B_t$ is ``built up'' from the
fibers corresponding to the generators appearing in the reduced
decomposition of $t$.  More precisely, if
  $t = x_1 x_2 \cdots x_n$ is in reduced form, then one requires that
  $B_t= B_{x_1} B_{x_2} \cdots B_{x_n}$ (meaning closed linear span).
This property makes sense for any group $G$, which, like the free group,
is equipped with a length function $|\cdot|$.  A Fell bundle over such a
group is said to be \ss if $B_{ts} = B_t B_s$ (closed linear span),
whenever $t$ and $s$ satisfy $|ts| = |t| + |s|$.

  Our main result, Theorem \lcite{6.3}, states, precisely, that any Fell
bundle over $\F$, which is orthogonal, \ss, and has separable fibers,
must be amenable.

To arrive at this conclusion we first restrict ourselves to a very
special case of Fell bundles, namely those which are associated to a
partial representation of $\F$ (see below for definitions).  For these
bundles, we prove an even stronger result, which is that they satisfy
the \ap of \cite{\Amena}.  This property implies amenability and also
some other interesting facts related to induced ideals of the
cross-sectional \cstar-algebra (see \scite{\Amena}{4.10}).

  The proof of the \ap for these restricted bundles is a direct
generalization of
  \scite{\Amena}{6.6}, where we proved that,
   for every \ss \prep $\pr$ of $\F_n$ (see below for definitions), such
that
  $
  \sum_{i=1}^n \pr(g_i)\pr(g_i)^* = 1,
  $
  the associated Fell bundle satisfies the approximation property. Here,
$\{g_1,\ldots,g_n\}$ are the generators of the free group $\F_n$.

Our generalization of this result, namely Theorem \lcite{3.7}, below,
amounts to replacing the hypothesis that $\sum_{i=1}^n
\pr(g_i)\pr(g_i)^* = 1$, by the weaker requirement that this sum is no
larger than $1$, or, equivalently, that the $\pr(g_i)\pr(g_i)^*$ are
pairwise orthogonal projections.

  Arriving at this generalization turns out to require a considerable
understanding of the various idempotents accompanying a \prep of $\F$,
and underlines the richness of ideas surrounding the concept of \preps.
In addition, the new hypothesis, that is, the orthogonality of these
projections, is easily generalizable to free groups with infinitely many
generators.  With the same ease, based on a simple inductive limit
argument, we extend Theorem \lcite{3.7} to the infinitely generated
case, obtaining Theorem \lcite{4.1}, below.

  In section 5, armed with this partial result, we study Fell bundles
which satisfy, in addition to the hypothesis of our main theorem, a
stability property.  Employing a fundamental result of Brown, Green and
Rieffel \cite{\BGR}, we are able to show that, for such bundles, there
is a hidden partial representation of $\F$ which sends us back to the
previously studied situation.  We finally remove the extra stability
hypothesis by means of a simple stabilization argument.

  It does not seem outlandish to expect that all amenable Fell bundles
satisfy some form of the \ap.  However, having no definite evidence that
this is so, we must be cautious in distinguishing these properties.
Accordingly, we must stress that our main result falls short of proving
the \ap for the most general situation treated, that is, of orthogonal
\ss bundles.  In this case, all we obtain is amenability, leaving that
stronger property as an open question.

  The author would like to express his thanks to the members of the
Mathematics Department at the University of Newcastle -- Australia,
where part of this work was developed.  Stimulating conversations with
M.~Laca and I.~Raeburn were essential to the formulation of the question
we have studied here.

  \section{Preliminaries}
  For the reader's convenience, and also to fix our notation, we shall
begin by briefly discussing some basic facts about partial group
representations, Fell bundles and the rich way in which these concepts
are interrelated.  The reader is referred to
  \cite{\FD},
  \cite{\Inverse},
  and
  \cite{\Amena}
  for more information on these subjects.

  Let $G$ be a group, fixed throughout this section.  Also, let $H$ be a
Hilbert space and denote by $\Lin(H)$ the algebra of all bounded linear
operators on $H$.

  \definition
  A \stress{\prep{}} of $G$ on $H$ is, by definition, a map
  $\pr : G \rightarrow \Lin(H)$
  such that
  \zitemno = 0
  \zitem $\pr(t)\pr(s)\pr(s\inv) = \pr(ts)\pr(s\inv)$,
  \zitem $\pr(t\inv) = \pr(t)^*$,
  \zitem $\pr(e) = I$,
  \medskip
  \noindent for all $t,s\in G$, where $e$ denotes the unit group element
and $I$ is the identity operator on $H$.

  Let $\pr$ be a \prep of $G$ on $H$.  It is an easy consequence of the
definition that each $\pr(t)$ is a partially isometric operator and
hence that
  $$
  \e(t) := \pr(t)\pr(t)^*
  $$
  is a projection (that is, a self-adjoint idempotent).  It is not hard
to show (see \cite{\Inverse}) that these projections commute among
themselves, and satisfy the commutation relation
  $$
  \pr(t)\e(s) = \e(ts)\pr(t),
  \eqno{(\seqnumbering)}\label{\CommRel}
  $$
  for all $t,s\in G$.

  There is a special kind of \preps worth considering, whenever $G$ is
equipped with a ``length'' function, that is, a non-negative real valued
function
  $
  |\cdot| : G \rightarrow {\bf R}_+
  $
  satisfying $|e| = 0$ and the triangular inequality
  $|ts| \leq |t|+|s|$.

  \definition
  A \prep $\pr$ of $G$ is said to be \stress{\ss{}} (with respect to a
given length function $|\cdot|$ on $G$) if
  $
  \pr(t)\pr(s) = \pr(ts)
  $
  whenever $t$ and $s$ satisfy
  $|ts| = |t| + |s|$.

  The concept of \preps is closely related to that of Fell bundles (also
known as \cstar-algebraic bundles \cite{\FD}) as we shall now see.

  \definition \label{\TheBundle}
  Given a \prep $\pr$ of $G$, for each $t$ in $G$, let
  $B^\pr_t$
  be the closed linear subspace of $\Lin(H)$ spanned by the set of
operators of the form
  $$
  \e(r_1) \e(r_2) \cdots \e(r_k) \pr(t),
  $$
  where $k\in\N$, and $r_1,r_2,\ldots,r_k$ are arbitrary elements of
$G$.

  Using the axioms of \preps and \lcite{\CommRel},
  it is an easy exercise to show (see \scite{\Amena}{Section 6}) that,
for all $t,s\in G$, we have
  $$
  B^\pr_t B^\pr_s \subseteq B^\pr_{ts}
  \and
  (B^\pr_t)^* = B^\pr_{t\inv}.
  $$
  Therefore, the collection
  $$
  \B^\pr := \{B^\pr_t\}_{t\in G}
  $$
  is seen to form a Fell bundle over $G$.

  \definition \label{\FellBundle}
  A Fell bundle over a discrete group $G$ is a collection
  $\B = \{B_t\}_{t\in G}$ of closed subspaces of $\Lin(H)$, such that
  $B_t B_s\subseteq B_{ts}$ and
  $(B_t)^* = B_{t\inv}$, for all $t$ and $s$ in $G$.

  As is the case with \cstar-algebras, which can be defined, concretely,
as a norm-closed *-subalgebra of $\Lin(H)$, as well as a certain
abstract mathematical object, defined via a set of axioms, Fell bundles
may also be seen under a dual point of view, specially if one restricts
attention to the case of discrete groups.  The above definition of Fell
bundles is the one we adopt here, referring the reader to
  \scite{\FD}{VIII.16.2} for the abstract version and to
  \scite{\FD}{VIII.16.4} for the equivalence of these.
  Nevertheless, it should be said that the point of view one usually
adopts in the study of Fell bundles stresses that each $B_t$ should be
viewed as a Banach space in its own, and that for each $t$ and $s$ in
$G$, one has certain algebraic operations
  $$
  \cdot : B_t \times B_s \rightarrow B_{ts}
  $$
  and
  $$
  * : B_t \rightarrow B_{t\inv},
  $$
  which, in our case, are induced by the multiplication and involution
on $\Lin(H)$, respectively.  If $G$ is not discrete, then one should
also take into account a topology on the disjoint union
  ${\buildrel \cdot \over \bigcup}_{t\in G}B_t$, which is compatible
with the other ingredients present in the situation.  See \cite{\FD} for
details.  Since we will only deal with Fell bundles over discrete
groups, we need not worry about this topology.

  \definition
  A Fell bundle $\B$, over of $G$, is said to be \stress{\ss{}} (with
respect to a given length function $|\cdot|$ on $G$) if
  $B_{ts} = B_t B_s$ (closed linear span), whenever $t$ and $s$ satisfy
  $|ts| = |t| + |s|$.

  When $\pr$ is \ss, one can prove, as in \scite{\Amena}{6.2}, that the
Fell bundle
  $\B^\pr$ is also \ss.

  \definition
  Given any Fell bundle
  $\B = \{B_t\}_{t\in G}$,
  with $G$ discrete,
  one defines its $l_1$ cross-sectional algebra
  \scite{\FD}{VIII.5},
  denoted $l_1(\B)$, to be the Banach *-algebra consisting of the $l_1$
cross-sections of $\B$, under the multiplication
  $$
  f g (t) = \sum_{s\in G} f(s) g(s\inv t) \for t\in G,\ f,g\in l_1(\B),
  $$
  involution
  $$
  f^*(t) = (f(t\inv))^* \for t\in G,\ f\in l_1(\B),
  $$
  and norm
  $$
  \|f\| = \sum_{s\in G} \|f(s)\| \for f\in l_1(\B).
  $$
  The cross-sectional \cstar-algebra of $\B$
  \scite{\FD}{VIII.17.2},
  denoted $\Cs$, is defined to be the enveloping \cstar-algebra of
$l_1(\B)$.

There is also a \stress{reduced} cross-sectional \cstar-algebra,
indicated by $\Csr$, which is defined to be the closure of $l_1(\B)$ in
a certain \stress{regular representation} (acting on the
right--$B_e$--Hilbert--bimodule formed by the $l_2$ cross-sections).
See
  \scite{\Amena}{2.3} for a precise definition.

  Both
  $\Cs$
  and
  $\Csr$
  contain a copy of the algebraic direct sum
  $\bigoplus_{t\in G} B_t$,
  as a dense subalgebra, making them into $G$-graded \cstar-algebras in
the sense of
  \scite{\FD}{VIII.16.11}
  (see also \scite{\Amena}{3.1}).  In both cases, the projections onto
the factors extend to bounded linear maps on the whole algebra, and, in
particular, for $B_e$, that projection gives a conditional expectation
  \scite{\Amena}{3.3}.

  This conditional expectation, say $E$, is faithful in the case of
  $\Csr$,
  in the sense that
  $$
  E(x^*x)=0\ \Rightarrow\ x=0
  $$
  for every $x\in \Csr$
  \scite{\Amena}{2.12}.
  However, the same cannot be said with respect to
  $\Cs$.
  In fact, there always exists an epimorphism
  $$
  \Lambda: \Cs \rightarrow \Csr
  $$
  which restrict to the identity map on
  $\bigoplus_{t\in G} B_t$ (see the discussion following
  \scite{\Amena}{2.2}
  as well as \scite{\Amena}{3.3}).
  The kernel of $\Lambda$ coincides with the \stress{degeneracy} ideal
for $E$, namely
  $$
  {\cal D} =
  \{ x\in \Cs : E(x^*x) = 0 \},
  $$
  where we are also denoting the conditional expectation for $\Cs$ by
$E$, by abuse of language
  (see \scite{\Amena}{3.6}).  This can be used to give an alternate
definition of $\Csr$, as the quotient of $\Cs$ by that ideal.

  The crucial property of Fell bundles with which we will be concerned
throughout this work is that of amenability.  This property is inspired,
first of all, in the corresponding concept for groups \cite{\Greenleaf},
but also in the work of Anantharaman-Delaroche \cite{\Claire} and Nica
\cite{\Nica}.  In the context we are interested, that is, for Fell
bundles, it first appeared in \cite{\Amena}, for discrete groups, and
was subsequently generalized by Ng for the non-discrete case \cite{\Ng}.

  \definition
  A Fell bundle $\B$, over a discrete group $G$, is said to be amenable,
if $\Lambda$ is an isomorphism.

  According to the characterization of the kernel of $\Lambda$, as in
our discussion above, we see that:

  \state Proposition \label{\FaithfullCondExp}
  A necessary and sufficient condition for $\B$ to be amenable is that
the conditional expectation $E$ of $\Cs$ be faithful.

  Any bundle is amenable when the base group $G$ is amenable
  \scite{\Amena}{4.7},
  whereas the typical example of non-amenable bundle is the
\stress{group bundle} \scite{\FD}{VIII.2.7} over a non-amenable group
$G$.  That is, the bundle $\C \times G$, with the operations
(abstractly) defined by
  $$
  (z,t) (w,s) = (zw,ts)
  \and
  (z,t)^* = (\bar z,t\inv)
  $$
  for $t,s\in G$ and $z,w\in \C$.  In this case $\Cs$ is the full group
\cstar-algebra of $G$, while $\Csr$ is its reduced algebra.  It is well
known
  \cite{\Ped}
  that $\Lambda$ is an isomorphism, in this case, if and only if $G$ is
amenable.

A Fell bundle may be amenable even if its base group $G$ is not.  One
example of this situation is given by
  \scite{\Amena}{6.7}.
  It consists of a Fell bundle over the non-amenable free group which
is, itself, amenable.  This example is particularly interesting, since
its cross-sectional \cstar-algebra is isomorphic to the Cuntz--Krieger
algebra ${\cal O}_A$.

  \definition \label{\APDef}
  We say that $\B$ has the \stress{\ap{}}
  \scite{\Amena}{4.5} if there exists a net
  $\{ a_i \}_{i\in I}$
  of finitely supported functions
  $a_i: G \arw B_e$,
  which is uniformly bounded in the sense that there exists a constant
$M>0$ such that
  $$
  \[ \sum_{t\in G} a_i(t)^* a_i(t) \] \leq M,
  $$
  for all $i$,
  and such that for all $b_t$ in each $B_t$ one has that
  $$
  b_t =
  \lim_{i \rightarrow \infty} \sum_{r\in G} a_i(tr)^* b_t a_i(r).
  $$

  The relevance of the \ap is that:

  \state Theorem \label{\APImpliesAmena}
  If a Fell bundle $\B$ has the approximation property, then it is
amenable.

  \proof See \scite{\Amena}{4.6}.
  \proofend

  For later use, it will be convenient to have certain equivalent forms
of the \ap, which we now study.

  \state Lemma \label{\Boundedness}
  Let $\B$ be a Fell bundle over $G$, and let
  $a: G \arw B_e$ be a finitely supported function.  Then, for each $t$
in $G$, the map
  $$
  b_t \in B_t \mapsto \sum_{r\in G} a(tr)^* b_t a(r) \in B_t
  $$
  is bounded, with norm no bigger than
  $\[ \sum_{r\in G} a(r)^* a(r) \]$.

  \proof Recall that
  $
  \[ \sum_{i=1}^n x_i^*y_i \] \leq
  \[ \sum_{i=1}^n x_i^*x_i \]\root
  \[ \sum_{i=1}^n y_i^*y_i \]\root,
  $
  whenever
  $x_1,\ldots,x_n$ and
  $y_1,\ldots,y_n$ are elements of a \cstar-algebra.  Therefore,
  letting $M = \[ \sum_{r\in G} a(r)^* a(r) \]$, we have,
  for all $b_t$ in $B_t$, that
  $$
  \[ \sum_{r\in G} a(tr)^* b_t a(r) \] \leq
  \[ \sum_{r\in G} a(tr)^* a(tr) \]\root
  \[ \sum_{r\in G} a(r)^* b_t^* b_t a(r) \]\root \$\leq
  M\root
  \|b_t\| \[ \sum_{r\in G} a(r)^* a(r) \]\root =
  M \|b_t\|.
  \proofend
  $$

  \state Proposition \label{\AlternateAP}
  Let $\B = \{B_t\}_{t\in G}$ be a Fell bundle over the discrete group
$G$.  Also, suppose we are given a dense subset $D_t$ of $B_t$, for each
$t$ in $G$.  Then the following are equivalent:
  \zitemno=0
  \zitem $\B$ satisfies the \ap.
  \zitem There exists a net
  $\{ a_i \}_{i\in I}$ satisfying all of the properties of
\lcite{\APDef}, except that the condition involving the limit is only
assumed for $b_t$ belonging to $D_t$.
  \zitem There exists a constant $M > 0$ such that, for all finite sets
  $
  \{b_{t_1},b_{t_2}, \ldots, b_{t_n} \},
  $
  with $b_{t_k}\in D_{t_k}$,
  and any $\varepsilon > 0$, there exists a finitely supported function
  $a: G \arw B_e$
  such that
  $
  \[ \sum_{t\in G} a(t)^* a(t) \] \leq M
  $ and $
  \[ b_{t_k} - \sum_{r\in G} a(t_k r)^* b_{t_k} a(r) \] < \varepsilon,
  $
  for $k=1,\ldots,n.$

  \proof The implications
  (i) $\Rightarrow$ (ii) and
  (ii) $\Rightarrow$ (iii) are obvious.  We shall than prove that
  (iii) $\Rightarrow$ (ii) and
  (ii) $\Rightarrow$ (i).
  With respect to our first task, consider the set of pairs
  $(X,\varepsilon)$,
  where $X$ is any finite subset of the disjoint union
  ${\buildrel \cdot \over \bigcup}_{t\in G}D_t$,
  and $\varepsilon$ is a positive real.  If these pairs are ordered by
saying that
  $(X_1,\varepsilon_1) \leq (X_2,\varepsilon_2)$
  if and only if
  $X_1 \subseteq X_2$ and $\varepsilon_1 \geq \varepsilon_2$,
  we clearly get a directed set.  For each such
  $(X,\varepsilon)$, let
  $a_{X,\varepsilon}$ be chosen such as to satisfy the conditions of
(iii) with respect to the set $X$ and $\varepsilon$.  It is then clear
that the net
  $\{ a_{X,\varepsilon}\}_{(X,\varepsilon)}$
  provides the required net.

  As for
  (ii) $\Rightarrow$ (i), since the maps of \lcite{\Boundedness} are
bounded, the convergence referred to in (ii) is easily seen to hold
throughout $B_t$.
  \proofend

  \section{Free groups and orthogonal \preps{}}
  This section is devoted to introducing a certain class of \preps of
the free group.  Let $\F$ denote the free group on a possibly infinite
set $\Gen$ (whose elements we call the generators).  Each $t$ in $\F$
has a unique decomposition (called its reduced decomposition, or reduced
form)
  $$
  t = x_1 x_2 \cdots x_k,
  $$
  where $x_i\in \Gen \cup \Gen\inv$ and $x_{i+1} \neq x_i\inv$ for all
$i$.
  In this case, we set $|t|=k$ and it is not hard to see that this
gives, in fact, a length function for $\F$.  It is with respect to this
length function that we will speak of \ss \preps of $\F$.

  \definition \label{\Ortho}
  A \prep $\pr$ of $\F$ is said to be \stress{orthogonal} if
  $\pr(x)^*\pr(y)=0$
  whenever $x,y\in \Gen$ are generators with $x\neq y$.

  A \prep of a group may not be determined by its values on a set of
generators.  For example, if we set
  $$
  \pr(t) = \left\{\matrix{
             1 & \hbox{if $|t|$ is even}\hfil\cr
             0 & \hbox{if $|t|$ is odd}\hfil\cr
           }\right.
  $$
  then $\pr$ is a \prep of $\F$ (on a one dimensional Hilbert space),
which coincides, on the generators, with the \prep{}
  $$
  \pr'(t) = \left\{\matrix{
             1 & \hbox{if $t=e$}\hfil\cr
             0 & \hbox{otherwise}\hfil\cr
           }\right.
  $$
  
  However, if $\pr$ is \ss, then
  $$
  \pr(t) = \pr(x_1) \pr(x_2) \cdots \pr(x_k),
  $$
  whenever $t = x_1 x_2 \cdots x_k$ is in reduced form.  Therefore, the
values of $\pr$ on the generators end up characterizing $\pr$
completely.  For this reason the fact that $\pr$ is orthogonal often
says little, unless one supposes that $\pr$ is \ss as well.

We shall denote by $\Pos$, the sub-semigroup of $\F$ generated by
$\Gen$, that is, the set of all products of elements from $\Gen$ (as
opposed to $\Gen\cup \Gen\inv$).  By convention, $\Pos$ also includes
the identity group element.  The elements of $\Pos$ are called the
\stress{positive} elements and will usually be denoted by letters taken
from the beginning of the Greek alphabet.  For each natural number $k$
we will denote by $\Pos_k$, the set of positive elements of length $k$.

  Note that, if $\pr$ is a \ss \prep of $\F$, and $\a,\b\in\Pos$, then
  $\pr(\a)\pr(\b)=\pr(\a\b)$, since
  $|\a\b|=|\a| + |\b|$.  This property will be useful in many
situations, below.

  \state Proposition \label{\PosNeg}
  Let $\pr$ be an orthogonal, \ss \prep of $\F$.  Then $\pr(t) = 0$ for
all elements $t$ in $\F$ which are not of the form $\mu\nu\inv$, with
$\mu$ and $\nu$ positive.

  \proof
  Let
  $t = x_1 x_2 \cdots x_k$,
  with $x_i\in \Gen \cup \Gen\inv$,
  be in reduced form.  Then, since $\pr$ is \ss,
  $\pr(t) = \pr(x_1) \pr(x_2) \cdots \pr(x_k)$.
  Now, because $\pr$ is orthogonal, if
  $x_i\in\Gen\inv$ and
  $x_{i+1}\in\Gen$ then
  $\pr(x_i) \pr(x_{i+1}) = 0$.
  So, in order to have $\pr(t)$ nonzero, all elements from $\Gen$ must
be to the left of the elements from $\Gen\inv$ in the decomposition of
$t$.  That is, $t$ is of the form described in the statement.
  \proofend

  \state Proposition \label{\POrth}
  Let $\pr$ be an orthogonal, \ss \prep of $\F$, and let $\a, \b \in
\Pos$.  If $|\a|=|\b|$, but $\a \neq \b$, then $\pr(\a)^*\pr(\b) = 0$.

  \proof Let $m=|\a|=|\b|$.  If $m=1$ then $\a$ and $\b$ are in $\Gen$
and the conclusion is a consequence of the orthogonality assumption.  If
$m > 1$ write
  $\a = x \tilde \a$ and
  $\b = y \tilde \b$ with
  $\tilde \a, \tilde \b \in \Pos$ and
  $x,y \in \Gen$.

  Assume, by way of contradiction, that
  $\pr(\a)^*\pr(\b) \neq 0$. Then
  $$
  0 \neq
  \pr(x \tilde \a)^* \pr(y \tilde \b) =
  \pr(\tilde \a)^* \pr(x)^* \pr(y) \pr(\tilde \b).
  $$
  So, in particular,
  $\pr(x)^* \pr(y) \neq 0,$
  which implies that $x=y$.

  We therefore have, using \lcite{\CommRel},
  $$
  0 \neq
  \pr(\tilde \a)^* \pr(x)^* \pr(x) \pr(\tilde \b) =
  \pr(\tilde \a)^* \e(x\inv) \pr(\tilde \b) =
  \pr(\tilde \a)^* \pr(\tilde \b) \e(\tilde\b\inv x\inv),
  $$
  which implies that $\pr(\tilde \a)^* \pr(\tilde \b) \neq 0$, and
hence, by induction, that $\tilde\a=\tilde\b$.  So $\a=\b$.
  \proofend
  
  The concept of orthogonality also applies to Fell bundles over free
groups.

  \definition
  A Fell bundle
  $\B = \{B_t\}_{t\in \F}$
  over $\F$ is said to be \stress{orthogonal} if
  $B_x^*B_y=\{0\}$
  whenever $x,y\in \Gen$ are generators with $x\neq y$.

  The parallel between this concept and its homonym \lcite{\Ortho} is
illustrated by our next:
  
  \state Proposition
  If $\pr$ is an orthogonal \prep of $\F$, then
  $\B^\pr$ is an orthogonal Fell bundle.

  \proof Left to the reader.
  \proofend

  \section{The finitely generated case}
  We now start the main technical section of the present work.  Here we
shall prove the \ap for Fell bundles arising from certain \preps of free
groups.  Even though our long range objective is to treat arbitrary free
groups, we shall temporarily restrict our attention to finitely
generated free groups.  So we make the following:

  \sysstate{Standing hypothesis}{\rm}{\label{\Standing}
  For the duration of this section,
  the set $\Gen$, of generators of $\F$ will be assumed to be finite and
$\pr$ will be a fixed
  orthogonal, \ss \prep of $\F$.}

Recall that $\e(t)$ denotes the final projection $\pr(t)\pr(t)^*$ of the
partial isometry $\pr(t)$.  In addition to $\e(t)$, the following
operators will play a crucial role:

  \def\newop#1#2{\bigskip$\displaystyle #1$, \quad $#2$}
  \newop{P_k = \sum_{\a\in\Pos_k} \e(\a)}{k\geq 1}
  \newop{Q_0 = 1 - P_1}{\null}
  \newop{f(t) = \pr(t) Q_0 \pr(t)^*}{t\in \F}
  \newop{Q_k = \sum_{\a\in\Pos_k} f(\a)}{k \geq 1}
  \bigskip

  The only place where the finiteness hypothesis in \lcite{\Standing}
will be explicitly used is in the observation that these sums are finite
sums.

  \state Proposition \label{\ManyProj}
  The following relations hold among the operators defined above:
  \zitemno=0
  \zitem $P_1$ and $Q_0$ are projections.
  \zitem Each $f(t)$ is a projection.
  \zitem If $t,s\in \F$ then $\pr(t)f(s)=f(ts)\pr(t)$.
  \zitem For every $t$ one has $f(t)\leq \e(t)$.
  \zitem If $\a,\b\in \Pos$ are such that $|\a|=|\b|$ but $\a \neq \b$,
then
  $\e(\a)\perp \e(\b)$,
  $f(\a)\perp f(\b)$,
  and
  $\e(\a)\perp f(\b)$.
  \zitem For all $k\geq 1$, both $P_k$ and $Q_k$ are projections and
  $Q_k = P_k - P_{k+1}$.
  \zitem For every $n$, we have that $Q_0 + Q_1 + \cdots + Q_{n-1} + P_n
= 1$.
  \zitem If $\a$ and $\b$ are distinct positive elements of $\F$ then,
regardless of their length, we have that $f(\a)\perp f(\b)$.

  \proof
  For $x,y\in\Pos_1=\Gen$, with $x\neq y$, we have, by the orthogonality
assumption, that
  $e(x)e(y) = \pr(x)\pr(x)^* \pr(y)\pr(y)^* = 0$.  Hence $P_1$ is a sum
of pairwise orthogonal projections, and thus, itself a projection.
Therefore $Q_0$ is also a projection.

  Speaking of (ii) we have
  $$
  f(t)^2 =
  \pr(t) Q_0 \pr(t)^* \pr(t) Q_0 \pr(t)^* =
  \pr(t) Q_0 \e(t\inv) Q_0 \pr(t)^*.
  $$
  Taking into account that the final and initial projections associated
to the partial isometries in a \prep all commute with each other
\cite{\Inverse}, we see that the above equals
  $$
  \pr(t) Q_0 \pr(t)^* \pr(t)\pr(t)^* =
  \pr(t) Q_0 \pr(t)^* =
  f(t).
  $$

  To prove (iii), let $t,s\in \F$.  Then
  $$
  \pr(t) f(s) =
  \pr(t) \pr(s) Q_0 \pr(s)^* =
  \pr(t) \pr(s) \pr(s)^* \pr(s) Q_0 \pr(s)^* \$=
  \pr(ts) \e(s\inv) Q_0 \pr(s)^* =
  \pr(ts) \pr(ts)^* \pr(ts) Q_0 \e(s\inv) \pr(s)^* \$=
  \pr(ts) Q_0 \pr(ts)^* \pr(ts) \pr(s)^* =
  \pr(ts) Q_0 \pr(ts)^* \pr(t) =
  f(ts)\pr(t).
  $$

  As for (iv)
  $$
  \e(t) f(t) =
  \pr(t) \pr(t)^* \pr(t) Q_0 \pr(t)^* =
  f(t).
  $$

  Given $\a$ and $\b$ as in (v) we have, by \lcite{\POrth}, that
  $
  \e(\a) \e(\b) =
  \pr(\a) \pr(\a)^* \pr(\b) \pr(\b)^* =
  0
  $
  which, when combined with (iv) above, yields the other statements of
(v).

  That $P_k$ is a projection follows from the fact that the summands in
its definition are pairwise orthogonal projections.  The same reasoning
applies to $Q_k$.  Now
  $$
  Q_k =
  \sum_{\a\in\Pos_k} \pr(\a) (1-P_1) \pr(\a)^* \$=
  \sum_{\a\in\Pos_k} \pr(\a) \pr(\a)^* -
    \sum_{\a\in\Pos_k} \sum_{x \in \Gen}
      \pr(\a) \pr(x) \pr(x)^* \pr(\a)^* \$=
  P_k - \sum_{\a\in\Pos_k} \sum_{x \in \Gen} \pr(\a x) \pr(\a x)^* \$=
  P_k - \sum_{\b\in\Pos_{k+1}} \pr(\b) \pr(\b)^* =
  P_k - P_{k+1}.
  $$

  To prove (vii), we just note that
  $$
  Q_0 + Q_1 + \cdots + Q_{n-1} + P_n \$=
  1 - P_1 \ + \ P_1 - P_2 + \cdots + P_{n-1} - P_n \ + \ P_n =
  1.
  $$

  Finally, let $\a\neq\b$ be positive and let $k=|\a|$ and $l=|\b|$.  If
$k=l$ then we have already seen in (v), that $f(\a)\perp f(\b)$.  On the
other hand, if $k\neq l$, then
  $$
  f(\a) \leq Q_k \perp Q_l \geq f(\b),
  $$
  where the orthogonality of $Q_k$ and $Q_l$ follows from (vii).  This
implies, again, that $f(\a)\perp f(\b)$.
  \proofend

  Recall that $\B^\pr = \{B^\pr_t\}_{t\in G}$ denotes the Fell bundle
associated to $\pr$, and consider, for each integer $n\geq 1$, the map
$b_n : \Pos \rightarrow B^\pr_e$ given by
  $$
  b_n(\a) = \left\{\matrix{
             f(\a) & \hbox{if $|\a| < n$}\hfil\cr
             \e(\a) & \hbox{if $|\a| = n$}\hfil\cr
             0 & \hbox{if $|\a| > n$}\hfil\cr
           }\right.
  $$

  \state Lemma \label{\SumB}
  For every $n\geq 1$ we have
  \quad $\displaystyle \sum_{\a\in\Pos} b_n(\a) = 1$.

  \proof We have
  $$
  \sum_{\a\in\Pos} b_n(\a) =
  \sum_{k=0}^n \sum_{\a\in\Pos_k} b_n(\a) =
  \sum_{k=0}^{n-1} \sum_{\a\in\Pos_k} f(\a) +
    \sum_{\a\in\Pos_n} \e(\a) =
  \sum_{k=0}^{n-1} Q_k + P_n = 1.
  \proofend
  $$

  The last relevant definition is that of another sequence of maps
  $
  a_n : \Pos \rightarrow B^\pr_e,
  $
  this time given by
  $$
  a_n(\a) = \({1\over n} \sum_{k=1}^n b_k(\a) \)\root.
  $$
  Note that $a_n(\a)=0$ for $|\a|>n$.
  We shall also think of the $a_n$ as functions defined on the whole of
$\F$, by setting $a_n(t)=0$ when $t$ is not positive.

  \state Lemma \label{\SumA}
  For every $n\geq 1$ we have
  \quad $\displaystyle \sum_{t\in\F} a_n(t)^* a_n(t) = 1$.

  \proof As already observed, $a_n(t)$ vanishes unless $t$ is positive.
In addition, $a_n(\a)$ is self-adjoint, so we must compute
  $$
  \sum_{\a\in\Pos} a_n(\a)^2 =
  \sum_{\a\in\Pos} {1\over n} \sum_{k=1}^n b_k(\a) =
  {1\over n} \sum_{k=1}^n \sum_{\a\in\Pos} b_k(\a) =
  1,
  $$
  where we have used \lcite{\SumB} in order to conclude the last step
above.
  \proofend

The square root appearing in the definition of $a_n$ can be explicitly
computed if we note that, for $|\a| \leq n$, we have
  $$
  \sum_{k=1}^n b_k(\a) =
  \(\sum_{k=|\a|+1}^n f(\a)\) + \e(\a) =
  (n - |\a|) f(\a) + \e(\a) \$=
  (n - |\a| + 1) f(\a) +
  \pef{\a},
  $$
  and that the expression above consists of a linear combination of
orthogonal projections, namely $f(\a)$ and $\e(\a) - f(\a)$ (see
\lcite{\ManyProj.iv}).  It follows that
  $a_n(\a)$ is given, explicitly, by
  $$
  a_n(\a) =
  \coef{\a} f(\a) +
  \({1 \over n}\)\root\pef{\a}.
  \eqno{(\seqnumbering)}\label{\NewAn}
  $$

   The following is the main technical point in showing the \ap for
$\B^\pr$:

  \state Lemma \label{\MainLemma}
  For every $t$ in $\F$ we have
  \quad $\displaystyle
  \pr(t) = \limn
           \sum_{r\in\F}
           a_n(tr)^* \pr(t) a_n(r).
  $

  \proof By \lcite{\PosNeg} we may assume that $t=\mu\nu\inv$, where
$\mu$ and $\nu$ are in $\Pos$.  We may also suppose that
$|t|=|\mu|+|\nu|$, that is, no cancelation takes place when $\mu$ and
$\nu\inv$ are multiplied together.

  In addition, since $a_n(r) = 0$ unless $r$ is positive and $|r| \leq
n$, each sum above is actually a finite sum.  In fact, the nonzero
summands in it are among those for which both
  $r$ and $tr$ are positive of length no bigger than $n$.

  Since $tr=\mu\nu\inv r$,
  if both $r$ and $tr$ are to be positive,
  we must have $r=\nu\b$, for some $\b\in\Pos$, and then
  $|tr| = |\mu\nu\inv\nu\b| = |\mu| + |\b|$.

  Also, in order to have $|r|$ and $|tr|$ no larger than $n$, we will
need
  $|r| = |\nu| + |\b| \leq n$, as well as
  $|\mu| + |\b| \leq n$, which are equivalent to
  $|\b| \leq m$,
  where
  $$
  m = \hbox{min}\{ n - |\nu|, n - |\mu|\}.
  $$
  Summarizing, for every $n$, we have
  $$
  \sum_{r\in\F} a_n(tr)^* \pr(t) a_n(r) =
  \sum_{|\b|\leq m} a_n(\mu\b) \pr(t) a_n(\nu\b),
  $$
  where we have also taken into account that each $a_n(t)$ is
self-adjoint.

  Substituting the expression for $a_n$, obtained in \lcite{\NewAn}, in
the above sum, we conclude that each individual summand equals
  $$
  \left[
    \coef{\mu\b} f(\mu\b) + \ninv
    \pef{\mu\b}
    \right] \$\cdot
  \pr(t)
  \left[
    \coef{\nu\b} f(\nu\b) + \ninv
    \pef{\nu\b}
    \right] =
  $$
  \bigskip
  $$
  =
  \coef{\mu\b} \coef{\nu\b} f(\mu\b) \ \pr(t) \ f(\nu\b) \$+
  \coef{\mu\b} \ninv f(\mu\b) \ \pr(t) \ \pef{\nu\b} \$+
  \ninv \coef{\nu\b} \pef{\mu\b} \ \pr(t) \ f(\nu\b) \$+
  \ninv \ninv \pef{\mu\b} \ \pr(t) \ \pef{\nu\b}.
  $$

  Let us indicate the four summands after the last equal sign above by
(i), (ii), (iii), and (iv), in that order.  In regards to (i), note
that, employing \lcite{\ManyProj.iii}, we have
  $$
  f(\mu\b) \pr(t) f(\nu\b) =
  f(\mu\b) \pr(\mu) \pr(\nu\inv) f(\nu\b) =
  \pr(\mu) f(\b) f(\b) \pr(\nu\inv) =
  \pr(\mu) f(\b) \pr(\nu)^*.
  $$

  Referring to (ii) we have
  $$
  f(\mu\b) \pr(t) \pef{\nu\b} =
  f(\mu\b) \pef{\mu\b} \pr(t) =
  0,
  $$
  because of \lcite{\ManyProj.iv}.  Similarly one proves that (iii)
vanishes as well.  As for (iv)
  $$
  \pef{\mu\b} \pr(t) \pef{\nu\b} \$=
  \pef{\mu\b} \pr(\mu) \pr(\nu\inv) \pef{\nu\b} =
  \pr(\mu) \pef{\b} \pr(\nu)^*.
  $$
  So,
  $$
  a_n(\mu\b) \pr(t) a_n(\nu\b) \$=
  \pr(\mu) \left[
  \coef{\mu\b} \coef{\nu\b} f(\b) +
  {1 \over n} \pef{\b}
  \right] \pr(\nu)^*,
  $$
  and the conclusion will follow once we prove that the term between
brackets above, summed over $|\b|\leq m$, converges, in norm, to the
identity operator, as $n\rightarrow\infty$.  We now set to do precisely
this.  Speaking of the identity operator, recall from \lcite{\SumA} and
\lcite{\NewAn} that,
  $$
  1 =
  \sum_{t\in\F} a_m(t)^2 =
  \sum_{|\b|\leq m} {m-|\b|+1 \over m} f(\b) + {1 \over m}\pef{\b}.
  $$
  Using this expression for the identity operator, we must then prove
the vanishing of the following limit
  $$
  \limn \Big\| \sum_{|\b|\leq m}
  \coef{\mu\b} \coef{\nu\b} f(\b) \$+
  {1 \over n} \pef{\b}
  - {m-|\b|+1 \over m} f(\b) - {1 \over m}\pef{\b}
  \Big\| \$=
  \limn \Big\| \sum_{|\b|\leq m}
  \({1 \over n} (n -|\mu\b| + 1)\root
    (n -|\nu\b| + 1)\root -
    {m-|\b|+1 \over m}\) f(\b) \$+
    \({1\over n} - {1\over m}\) \pef{\b} \Big\| \$\leq
  \limn \[ \sum_{|\b|\leq m}
  \({1 \over n} (n -|\mu\b| + 1)\root
    (n -|\nu\b| + 1)\root -
    {m-|\b|+1 \over m}\) f(\b) \] \$+
  \limn \[ \sum_{|\b|\leq m}
  \({1\over n} - {1\over m}\) \pef{\b}
  \].
  $$

  The two limits will now be shown to equal zero.  We should point out
that, with respect to the the first one, we are facing a linear
combination of pairwise orthogonal projections by
\lcite{\ManyProj.viii}.  The same, however, is not true for the second.

  Using this observation, we see that the norm, in the first case,
equals
  $$
  \max_{|\b|\leq m} \left| {1 \over n}
  (n -|\mu\b| + 1)\root
  (n -|\nu\b| + 1)\root -
  {m -|\b| + 1\over m}
  \right|.
  $$
  In order to show that this goes to zero as $n \rightarrow \infty$, let
us assume, without loss of generality, that
  $|\mu| \geq |\nu|$, and hence that
  $m=n-|\mu|$.

In addition, it is easy to see that, for every pair of positive reals
$x$ and $y$, one has that
  $|x-y| \leq |x^2 - y^2|\root$.  So, the task facing us can be replaced
by
  $$
  \limn \max_{|\b|\leq m} \left|
  {(n -|\mu\b| + 1) (n -|\nu\b| + 1) \over n^2} -
  {(n -|\mu\b| + 1)^2 \over (n -|\mu|)^2}
  \right| \={?} 0.
  $$
  The term between the single bars is no bigger than
  $$
  \left| {(n - |\mu\b| + 1) (n - |\nu\b| + 1) \over n^2} -
    {(n - |\mu\b| + 1)^2 \over n^2} \right| \$+
  (n - |\mu\b| + 1)^2
  \left| {1 \over n^2} - {1 \over (n - |\mu|)^2} \right| \$\leq
  {n-|\mu\b| + 1 \over n^2} \Big| |\mu|-|\nu| \Big|
  +
  (n-|\mu\b| + 1)^2\left|{ -2n|\mu| + |\mu|^2 \over
    n^2(n-|\mu|)^2 }\right| \$\leq
  {n + 1 \over n^2} \Big| |\mu|-|\nu| \Big|
  +
  (n+1)^2\left|{ -2n|\mu| + |\mu|^2 \over
    n^2(n-|\mu|)^2 }\right|,
  $$
  which is now easily seen to go to zero, uniformly on $\b$, as
$n\rightarrow\infty$.

  To conclude, we need only show the vanishing of
  \def\pre{\limn \left|{1\over n} - {1\over m}\right|}
  $$
  \limn \[ \sum_{|\b|\leq m}
    \({1\over n} - {1\over m}\) \pef{\b} \] =
  \pre \[ \sum_{k=0}^m \sum_{\b\in\Pos_k}\e(\b) -f(\b)\] \$=
  \pre \[ \sum_{k=0}^m P_k -Q_k \] =
  \pre \[ \sum_{k=0}^m P_{k+1} \] \leq
  \pre (m+1).
  $$

 Now, recalling our assumption that $m=n-|\mu|$, the limit above equals
  $$
  \limn \left| {1\over n} -{1\over n-|\mu|} \right| (n-|\mu|+1) =
  \limn {|\mu|(n-|\mu|+1) \over n(n-|\mu|)} =
  0.
  \proofend
  $$

  We are now prepared to face one of our main goals.

  \state Theorem \label{\MainFinite}
  Let $\pr$ be an orthogonal, \ss \prep of a finitely generated free
group $\F$.  Then the Fell bundle $\B^\pr$ satisfies the approximation
property and hence is amenable.  Moreover, the constant $M$ referred to
in \lcite{\APDef} may be taken to be $1$.

  \proof Let $\{a_n\}_{n\in \N}$ be defined as above.  Then
  $
  \[ \sum_{t\in \F} a_n(t)^* a_n(t) \] = 1,
  $
  for all $n$, by \lcite{\SumA}, and employing \lcite{\AlternateAP.ii},
it is now enough to show that
  $$
  b_t =
  \lim_{n \rightarrow \infty} \sum_{r\in \F} a_n(tr)^* b_t a_n(r),
  $$
  for all $b_t$ of the form
  $b_t = \e(r_1) \e(r_2) \cdots \e(r_k) \pr(t)$,
  where $k\in\N$, and $r_1,r_2,\ldots,r_k$ are arbitrary elements of
$\F$.  This is because the linear combinations of the elements of this
form are dense in $B^\pr_t$, by definition \lcite{\TheBundle}.

We have already observed that the projections associated to the partial
isometries in a \prep form a commutative set.  Since $a_n(\a)$ is given
by a linear combination of such projections, by \lcite{\NewAn} (and
$a_n(t)=0$ when $t\notin \Pos$), it is clear that $a_n(t)$ commutes with
the $\e(r_j)$.  Therefore, by \lcite{\MainLemma},
  $$
  \lim_{n \rightarrow \infty} \sum_{r\in \F} a_n(tr)^* b_t a_n(r) \$=
  \e(r_1) \e(r_2) \cdots \e(r_k)
  \lim_{n \rightarrow \infty} \sum_{r\in \F} a_n(tr)^* \pr(t) a_n(r) \$=
  \e(r_1) \e(r_2) \cdots \e(r_k) \pr(t) =
  b_t.
  $$
  This concludes the proof of the \ap and hence of the amenability of
$\B^\pr$, by \lcite{\APImpliesAmena}.
  \proofend

  \section{Arbitrary free groups}
  We will now extend the results of the previous section, by dropping
the finiteness hypothesis of \lcite{\Standing}, and hence including
infinitely generated free groups in our study. The strategy will be to
adapt the work done above, to the general case, using an inductive limit
argument.

  Let, therefore, $\F$ be the free group on a set $\Gen$, no longer
assumed to be finite, or even countable.  Also let $\pr$ be an
orthogonal, \ss \prep of $\F$, considered fixed throughout this section.

  For each finite subset $X$ of $\Gen$, let $\F_X$ denote the subgroup
of $\F$ generated by $X$.  It is quite obvious that $\F_X$ is again a
free group, and that $\F$ is the union of the increasing net
  $\{\F_X\}_X$.
 The length functions we've been considering are compatible in the sense
that the one for $\F$ restricts to the one for $\F_X$.  Therefore the
restriction of $\pr$ to $\F_X$ is also \ss, and obviously also
orthogonal.  Let $\B^X$ denote the Fell bundle for $\pr|_{\F_X}$, as in
\lcite{\TheBundle}.

It is clear that, for each $t$ in $\F$, one has that $B^\pr_t$ is the
closure of the union of the $B^X_t$, as $X$ ranges in the collection of
finite subsets $X\subseteq \Gen$, such that $t\in \F_X$.

  \state Theorem \label{\ArbitraryAmenable}
  Let $\pr$ be an orthogonal, \ss \prep of an arbitrary free group $\F$.
Then the Fell bundle $\B^\pr$ satisfies the approximation property and
hence is amenable.

  \proof For each $t$ in $\F$, let $D_t$ be the union of the
  $B^X_t$, as described above, which is dense in $B^\pr_t$.  We will now
prove \lcite{\AlternateAP.iii}, with respect to this choice of $D_t$.
  Let $M=1$.  Then, given a finite set
  $ \{b_{t_1},b_{t_2}, \ldots, b_{t_n} \}, $
  with $b_{t_k}\in D_{t_k}$,
  and any $\varepsilon > 0$, there clearly exists a single finite
$X\subseteq \Gen$, such that every
  $b_{t_k}\in B^X_{t_k}$.
  Now, by \lcite{\MainFinite} we conclude that a finitely supported map
  $
  a: \F_X \arw B^X_e \subseteq B^\pr_e
  $
  exists, satisfying
  $
  \[ \sum_{t\in \F_X} a(t)^* a(t) \] \leq 1,
  $
  and
  $
  \[ b_{t_k} - \sum_{r\in \F_X} a(t_k r)^* b_{t_k} a(r) \] <
\varepsilon,
  $
  for all $k=1,\ldots,n$.
  If we extend $a$ to the whole of $\F$ by declaring it zero outside
$\F_X$, then these two sums may be taken for ${t\in \F}$, as opposed to
${t\in \F_X}$, without changing the result.  We conclude that
\lcite{\AlternateAP.iii} holds, and hence that $\B^\pr$ satisfies the
\ap.
  \proofend

  \section{Stable Fell bundles}
  We shall now treat Fell bundles, over arbitrary free groups, which are
not necessarily presented in terms of a \prep.  This section does not
yet contain our strongest result, because we shall be working under the
assumption that the unit fiber algebra of $\B$, that is $B_e$, is a
stable \cstar-algebra, in the sense that it isomorphic to
$B_e\otimes\K$, where $\K$ is the algebra of compact operators on a
separable, infinite dimensional Hilbert space.  In the next and final
section we will then remove this stability hypothesis.

  The method we shall adopt here will be to construct a \prep which is
closely related to $\B$.  In fact this method is, essentially, the one
we have used in \cite{\TPA} to obtain the classification of Fell bundles
in terms of twisted partial actions (see \scite{\TPA}{7.3}).  However,
we will refrain from utilizing the full machinery of \cite{\TPA} for two
reasons.  First, by proceeding more or less from scratch, and using the
special features of the free group, we will be able to make the
presentation somewhat more elementary and self contained.  Secondly, the
classification theorem mentioned includes a 2-cocycle, which will not
show up here, again because of the special properties of the group we
are dealing with.

  We begin with some simple facts about partial isometries and
projections on a Hilbert space $H$.

  \state Lemma \label{\Idempot}
  Let $p$ be an operator on $H$, such that $p^2=p$ and $\|p\|\leq
1$. Then $p=p^*$.

  \proof Let $\xi\in p(H)^\perp$.  Then, for every $\lambda\in\R$ we
have
  $
  \|p(\xi+\lambda p(\xi))\| \leq \|\xi+\lambda p(\xi)\|,
  $
  which implies, after a short calculation, that
  $
  (1+2\lambda)\|p(\xi)\|^2 \leq \|\xi\|^2,
  $
  and hence that $p(\xi)=0$.  This says that $p$ vanishes on
$p(H)^\perp$ and, since $p$ is the identity on $p(H)$, then it must be
the orthogonal projection onto $p(H)$.  Hence $p=p^*$.
  \proofend

  \state Lemma \label{\ComutProj}
  Let $p$ and $q$ be projections (self-adjoint idempotents) in
$\Lin(H)$.  Then $pq$ is idempotent if and only if
  $p$ and $q$ commute.
  
  \proof If $pq$ is idempotent, then, since $\|pq\|\leq 1$, we have, by
\lcite{\Idempot}, that $pq = (pq)^* = qp$.  The converse is trivial.
  \proofend

  \state Lemma \label{\ProdPIso}
  Let $u$ and $v$ be partial isometries in $\Lin(H)$.  Then $uv$ is a
partial isometry if and only if $u^*u$ and $vv^*$ commute (see also
\cite{\Nandor}).

  \proof We have that $uv$ is a partial isometry, if and only if
  $$
  uv(uv)^*uv = uv \ \Longleftrightarrow\
  uvv^*u^*uv = uv \$\ \Longleftrightarrow\
  u^*uvv^*u^*uvv^* = u^*uvv^* \ \Longleftrightarrow\
  (u^*uvv^*)^2 = u^*uvv^*,
  $$
  which, by \lcite{\ComutProj}, is equivalent to the commutativity of
  $u^*u$ and $vv^*$.
  \proofend

  \state Proposition \label{\PisoSet}
  Let $U = \{u_x\}_{x\in\Gen}$ be a family of partial isometries on a
Hilbert space $H$ and let $\I$ be the multiplicative sub-semigroup of
$\Lin(H)$ generated by $U \cup U^*$. Denote by $\F$ the free group on
$\Gen$.  Then, the following are equivalent:
  \zitemno=0
  \zitem There exists a \ss \prep $\pr$ of $\F$ such that $\pr(x)=u_x$
for every $x\in\Gen$.
  \zitem There exists a \prep $\pr$ of $\F$ such that $\pr(x)=u_x$ for
every $x\in\Gen$.
  \zitem Every $u$ in $\I$ is a partial isometry.
  \zitem For any $u,v\in\I$ we have that $uu^*$ and $vv^*$ commute.

  \proof
  (i) $\Rightarrow$ (ii): Obvious.

  (ii) $\Rightarrow$ (iii):
  Recall that an operator $u$ is a partial isometry, if and only if
$uu^*u=u$.
  So, let $u = u_1\cdots u_n$, with $u_i\in U\cup U^*$, and take $t_i\in
\Gen\cup\Gen\inv$ such that $\pr(t_i)=u_i$.  Then
$u=\pr(t_i)\cdots\pr(t_n)$ and, by induction on $n$,
  $$
  uu^*u
  =
  \pr(t_1)\cdots\pr(t_{n-1})
  \e(t_n)
  \pr(t_{n-1})^*\cdots\pr(t_1)^*
  \pr(t_1)\cdots\pr(t_n)
  \$=
  \e(t_1\cdots t_n)
  \pr(t_1)\cdots\pr(t_{n-1})
  \pr(t_{n-1})^*\cdots\pr(t_1)^*
  \pr(t_1)\cdots\pr(t_n)
  \$=
  \e(t_1\cdots t_n)
  \pr(t_1)\cdots\pr(t_{n-1})
  \pr(t_n)
  \$=
  \pr(t_1)\cdots\pr(t_{n-1})
  \e(t_n)
  \pr(t_n)
  =
  \pr(t_1)\cdots\pr(t_{n-1})
  \pr(t_n)
  =
  u.
  $$

  (iii) $\Rightarrow$ (iv):
  If $u,v\in\I$, then $u^*v\in\I$ and, by (iii), it is a partial
isometry.  Hence, using \lcite{\ProdPIso}, we have that $uu^*$ and
$vv^*$ commute.

  (iv) $\Rightarrow$ (i): Define, for all $x \in \Gen$,
  $\pr(x)=u_x$ and
  $\pr(x\inv)=u_x^*$.  Now, if $t= x_1\cdots x_n$, with
  $x_i\in\Gen\cup\Gen\inv$, is in reduced form, put
  $\pr(t)=\pr(x_1)\cdots\pr(x_n)$.
  It is then obvious that
  $
  \pr(t)\pr(s) = \pr(ts)
  $
  whenever $t$ and $s$ satisfy
  $|ts| = |t| + |s|$.

We claim that $\pr$ is a \prep of $\F$.  The crucial point is to prove
that
  $
  \pr(t) \pr(s) \pr(s)^* = \pr(ts) \pr(s)^*
  $
  for $t,s \in \F$.
  To do this we use induction on
  $|t| + |s|$.  If either $|t|$ or $|s|$ is zero, there is nothing to
prove.  So, write
  $t = \tilde t x$ and
  $s = y \tilde s$ ,
  where $x,y\in \Gen \cup \Gen\inv$ and, moreover,
  $|t|=|\tilde t|+1$ and
  $|s|=|\tilde s|+1$.

  In case
  $x\inv \neq y$ we have
  $|ts|=|t| + |s|$ and hence $\pr(ts) = \pr(t) \pr(s)$.
  If, on the other hand,
  $x\inv = y$, we have
  $$
  \pr(t) \pr(s) \pr(s)^* =
  \pr(\tilde t x) \pr(x\inv \tilde s) \pr(\tilde s\inv x) =
  \pr(\tilde t) \pr(x) \pr(x)^* \pr(\tilde s) \pr(\tilde s)^* \pr(x).
  $$
  By (iv) and the induction hypothesis, we conclude that the above
equals
  $$
  \pr(\tilde t) \pr(\tilde s) \pr(\tilde s)^* \pr(x) \pr(x)^* \pr(x) =
  \pr(\tilde t \tilde s) \pr(\tilde s)^* \pr(x)
  \$=
  \pr(ts) \pr(\tilde s\inv x) =
  \pr(ts) \pr(s)^*.
  \proofend
  $$
 
  It would be interesting to find a condition about a set $U$ of partial
isometries, which is equivalent to the ones above, but which refers
exclusively to the $u_x$'s, themselves, rather than to arbitrary
products of them.  A related observation is that the Cuntz--Krieger
\cite{\CK} relations were shown to imply the conditions above
\scite{\Amena}{5.2}.

  In our next result, we will use two important concepts from
\cite{\TPA}, which we briefly summarize here.
  By definition, a TRO (for ternary ring of operators), is a closed
linear subspace $E\subseteq \Lin(H)$, such that $EE^*E\subseteq E$
  (see also \cite{\Zettl}).
  We adopt the convention that the product of two or more sets, as
above, is supposed to mean the \stress{closed linear span} of the set of
products.

Given a TRO $E$, we say that a partial isometry $u$ is associated to $E$
  \scite{\TPA}{5.4},
  and write
  $u \assoc E$,
  if $u^*E = E^*E$ and $uE^*=EE^*$.  If, in addition, the range of the
final projection $uu^*$ coincides with $EH$ (equivalently, if the range
of the initial projection $u^*u$ coincides with $E^*H$), then we say
that $u$ is \stress{strictly} associated to $E$
  \scite{\TPA}{5.5},
  and write
  $u \sassoc E$.

  It is a consequence of
  \scite{\BGR}{3.3 and 3.4}, that whenever $E$ is separable and stable,
in the sense that $E^*E$ and $EE^*$ are stable \cstar-algebras
  \scite{\TPA}{4.11}, then a partial isometry strictly associated to $E$
always exists (see also \scite{\TPA}{5.3 and 5.2}).

   The reason why TROs are relevant here is that any fiber of a Fell
bundle is a TRO, as one may easily see.  But, because of the
separability requirement, we shall now restrict to the case of separable
bundles, according to the following:

  \definition
  A Fell bundle $\B$, over a discrete group $G$, is said to be
\stress{separable} if each $B_t$ is a separable Banach space.

  \state Theorem \label{\ExistPrep}
  Let $\B = \{B_t\}_{t\in \F}$ be a \ss, separable Fell bundle over the
arbitrary free group $\F$, represented on a Hilbert space $H$.  Suppose
that $B_e$ is stable.  Then, there exists a \ss \prep{} $\pr$ of $\F$ on
$H$, such that $\pr(t)\sassoc B_t$.  In addition, if $\B$ is orthogonal,
then $\pr$ is necessarily orthogonal, as well.

  \proof By
  \scite{\TPA}{4.12}, each $B_t$ is stable and hence, by
  \scite{\TPA}{5.2}, there exists a partial isometry $u_t\sassoc B_t$.
Let
  $U = \{u_x\}_{x\in\Gen}$, where, as before, $\Gen$ is the set of
generators of $\F$.  We claim that $U$ satisfies \lcite{\PisoSet.iii}.
In fact, we shall prove that, given
  $x_1,\ldots,x_n$ in $\Gen$, then
  $B_{x_1}\cdots B_{x_n}$ is a TRO and that
  $u_{x_1}\cdots u_{x_n}$ is a partial isometry strictly associated to
it.  Proceeding by induction, we may assume that
  $E = B_{x_1}\cdots B_{x_{n-1}}$ is a TRO and that
  $u = u_{x_1}\cdots u_{x_{n-1}}\sassoc E$.  Now, observe that
  $E^*E$ and $B_{x_n}B_{x_n}^*$ are ideals in $B_e$, and hence they
commute, as sets, that is,
  $E^*EB_{x_n}B_{x_n}^* = B_{x_n}B_{x_n}^*E^*E$.  So, by
\scite{\TPA}{6.4}, it follows that $EB_{x_n}$ is a TRO, and that
  $uu_{x_n}\sassoc EB_{x_n}$.  This proves our claim.  So, let $\pr$ be
a \ss \prep of $\F$ satisfying \lcite{\PisoSet.i}.  It remains to show
that
  $\pr(t) \sassoc B_t$, but this follows from the conclusion just above,
once we write $t=x_1\ldots x_n$ in reduced form.
  
  Assume, now, that $\B$ is orthogonal.  Then, given $x\neq y$, in
$\Gen$, we have, again by \scite{\TPA}{6.4}, that $\pr(x)^*\pr(y)\sassoc
B_x^* B_y = \{0\}$.  Therefore $\pr(x)^*\pr(y) = 0$.
  \proofend

  \state Theorem \label{\StableIsAmena}
  Let $\B$ be an orthogonal, \ss, separable Fell bundle over $\F$, and
suppose that $B_e$ is stable.  Then $\B$ is amenable (see below for the
non-stable case).

  \proof Let $H$ be the space where $\B$ acts, and let $\pr$ be the
orthogonal, \ss \prep of $\F$ on $H$, provided by \lcite{\ExistPrep}.
As usual, we denote by $\e(t)$ the final projection of $\pr(t)$.

  Let $\B^\pr$ be the Fell bundle associated to $\pr$ as in
\lcite{\TheBundle}.  By \lcite{\ArbitraryAmenable}, we know that
$\B^\pr$ satisfies the \ap.  Let, therefore, $\{ a_i \}_{i\in I}$ be a
net of maps satisfying the conditions of \lcite{\APDef}, with respect to
$\B^\pr$.

  We claim that $\e(t)$ commutes with $B_e$.  In fact, because
$\pr(t)\sassoc B_t$, we have that $\e(t)$ is the orthogonal projection
onto $B_tH$.  Now, observing that $B_e$ leaves the latter invariant, we
obtain the conclusion.
  It follows that the \cstar-algebra generated by all the $\e(t)$,
namely $B_e^\pr$, is contained in the commutant of $B_e$.

  Let $t\in\F$ and pick $b_t\in B_t$.  Define
  $c_t = b_t\pr(t)^*$.  Then
  $
  c_t \in B_t\pr(t)^* = B_t B_t^* \subseteq B_e.
  $
  On the other hand, note that, since the range of $b_t^*$ is contained
in $B_{t\inv}H$, which is also the range of $\e(t\inv)$, we have that
$\e(t\inv)b_t^* = b_t^*$, or simply
  $b_t\e(t\inv)= b_t$.  This implies that
  $
  c_t\pr(t) =
  b_t\pr(t)^*\pr(t) =
  b_t\e(t\inv) =
  b_t.
  $

  Since each $a_i(r)\in B_e^\pr$, we have that it commutes with $c_t$
and hence
  $$
  \lim_{i \rightarrow \infty} \sum_{r\in \F} a_i(tr)^* b_t a_i(r) =
  c_t \lim_{i \rightarrow \infty} \sum_{r\in \F} a_i(tr)^* \pr(t) a_i(r)
=
  c_t \pr(t) =
  b_t.
  $$
  
  We therefore see that the net $\{a_i\}$ satisfies the properties of
\lcite{\APDef} with respect to $\B$, except that there is no reason to
expect that the values of the maps $a_i$ lie in $B_e$.  Therefore this
falls short of proving the \ap for $\B$ and hence we cannot use
  \lcite{\APImpliesAmena} to conclude the amenability of $\B$.

  Fortunately, what we do have is enough to fit the hypothesis of a
slight generalization of the results of \cite{\Amena} leading to
\lcite{\APImpliesAmena}, which goes as follows: Let $\Cs$ be faithfully
represented on a Hilbert space $K$.  Since $\Cs$ contains
$\bigoplus_{t\in G} B_t$, we may then identify each $B_t$ with its image
under that representation, and then think of $B_t$ as a space of
operators on $K$.  In other words, this provides a faithful
representation of $\B$ as operators on $K$, and hence we might as well
assume that $H=K$, which we do, from now on. Under this assumption, we
have that the sub-\cstar-algebra of $\Lin(H)$ generated by $\bigcup_t
B_t$ is isomorphic to $\Cs$.

  For each $t$, consider the space $C_t$ of operators on $H\*l_2(\F)$,
given by $C_t=B_t\*\lambda_t$, where $\lambda$ is the left regular
representation of $\F$.  It is easy to see that the $C_t$ form a Fell
bundle, which is again isomorphic to $\B$.  However, the \cstar-algebra
generated by $\bigcup_t C_t$ is now isomorphic to $\Csr$, a fact that
follows from
  \scite{\Amena}{3.7}.

  Recall that the reasoning at the beginning of the present proof
provided a net
  $\{ a_i \}_{i\in I}$
  of maps
  $a_i: G \arw \Lin(H)$, such that
  $
  \sup_i \[ \sum_{t\in G} a_i(t)^* a_i(t) \] < \infty,
  $
  and that
  $
  b_t = \lim_{i \rightarrow \infty} \sum_{r\in G} a_i(tr)^* b_t a_i(r),
  $
  for all $b_t$ in each $B_t$.  Following the argument used in
  \scite{\Amena}{4.3},
  let, for each $i$,
  $$
  V_i : H \arw H\*l_2(\F),
  $$
  be given by the formula
  $V_i(\xi) = \sum_{t\in \F} a_i(t) \xi\*\delta_t,$
  where $\{\delta_t\}$ is the standard orthonormal base of $l_2(\F)$.
One then proves, as in
  \scite{\Amena}{4.3}, that
  $\| V_i \| \leq \| \sum_{t\in \F} a_i(t)^* a_i(t) \|^{1/2}$, and hence
that the $V_i$ are uniformly bounded.

  Now, define the completely positive maps
  $$
  \Psi_i : \Lin(H\* l_2(\F)) \arw \Lin(H)
  $$
  by $\Psi_i(T) = V_i^* T V_i$, for each $T$ in $\Lin(H\* l_2(\F))$.
Again as in \scite{\Amena}{4.3}, one has that, for every $b_t$ in $B_t$,
  $$
  \Psi_i(b_t\*\lambda_t) = \sum_{r\in \F} a_i(tr)^* b_t a_i(r).
  $$

  The somewhat annoying fact that $a_i(t)$ may not belong to $B_e$
forbids us to say that $\Psi_i$ is a map from $\Csr$ into $\Cs$,
  as stated in
  \scite{\Amena}{4.3}.  This, however, is not a cause for despair.

  Consider the canonical map
  $
  \Lambda : \Cs \arw \Csr,
  $
  described in section 1.  Under the current representation of $\Csr$ on
$H\* l_2(\F)$, we have that $\Lambda(b_t) = b_t\*\lambda_t$ for all
$b_t$.
  Now, consider the composition of maps
  $$
  \Cs \ {\buildrel \Lambda \over \longrightarrow} \
  \Csr \ {\buildrel \Psi_i \over \longrightarrow} \
  \Lin(H).
  $$
  For $b_t$ in $B_t$, we have
  $$
  \lim_i \Psi_i(\Lambda(b_t)) =
  \lim_i \Psi_i(b_t\*\lambda_t) =
  \lim_i \sum_{r\in \F} a_i(tr)^* b_t a_i(r) =
  b_t.
  $$
  Now, by the uniform boundedness of these maps we then conclude that
  $\lim_i \Psi_i(\Lambda(x)) = x$
  for all $x\in\Cs$.  This implies that $\Lambda$ is injective and hence
that $\B$ is amenable, as required.
  \proofend

  \section{The general case}
  In this section we will prove our most general result, which is the
amenability of orthogonal, \ss, separable bundles.  This amounts to
dropping the stability hypothesis of the previous section, which we do
by a ``stabilization argument''.  Ideally one should develop the whole
theory of tensor products for Fell bundles but we feel this is not the
place to do it.  Instead, we perform our tensor products in a way which
is enough for our purposes, albeit in a somewhat crude manner.

  Let $\B$ be any Fell bundle over a discrete group $G$, acting on the
Hilbert space $H$.  As in the proof of
  \lcite{\StableIsAmena},
  we may assume that the sub-\cstar-algebra of $\Lin(H)$ generated by
  $\bigcup_t B_t$
  is isomorphic to $\Cs$.

For each $t\in G$, consider the the subset of
  $\Lin(H\*l_2)$
  (where $l_2$ is the usual infinite dimensional separable Hilbert
space), denoted by $B_t\*\K$, and defined by
  $$
  B_t\*\K = \overline{\hbox{span}} \{ b_t\otimes k : b_t\in B_t, k\in \K
\}.
  $$
  Here $\K$ is the algebra of compact operators on $l_2$.
  It is elementary to show that
  $\B\*\K$, as defined by
  $\B\*\K = \{B_t\*\K\}_{t\in G}$,
  is a Fell bundle in its own right.  We shall say that $\B\*\K$ is the
\stress{stabilization} of $\B$.

  \state Proposition \label{\Kiso}
  Let $\B$ be a Fell bundle.  Then
  $C^*(\B\*\K)$ is isomorphic to $\Cs\*\K$.

  \proof
  Let us temporarily use the notation
  $\tilde\B$ for $\B\*\K$ and $\tilde B_t$ for $B_t\*\K$.
  Observe that $\K$ may be viewed, in a canonical way, as a subalgebra
of the multiplier algebra $\M(\tilde B_e)$, which, in turn, may be
viewed as a subalgebra of $\M(C^*(\tilde\B))$
  \scite{\FD}{VIII.5.8 and 1.15}.
  Since one clearly has
  $
  \K C^*(\tilde\B) =
  C^*(\tilde\B),
  $
  one may now show that
  $C^*(\tilde\B)$ is isomorphic to $A\*\K$, where
  $
  A=p C^*(\tilde\B) p,
  $
  and $p$ is any minimal projection in $\K$.

  Using the universal property
  \scite{\FD}{VIII.16.12}
  of cross-sectional \cstar-algebras, one may show that the assignment
  $b_t\in B_t \mapsto b_t\*p \in C^*(\tilde\B)$
  extends to a surjective *-homomorphism
  $\phi : \Cs \rightarrow A$.

  We claim that $\phi$ is injective.  In fact,
  let $H$ be the space where $\B$ acts, so that $\tilde\B$ sits inside
of $\Lin(H\*l_2)$.  The universal property, this time applied to
$\tilde\B$, implies that there exists a *-representation
  $$
  \pi : C^*(\tilde\B) \rightarrow \Lin(H\*l_2)
  $$
  which, restricted to each $\tilde B_t$, coincides with the inclusion
of $\tilde B_t$ in $\Lin(H\*l_2)$.  It is then easy to see that
$\pi\phi$ maps $\Cs$ onto the closed linear span of $\bigcup_t B_t\*p$,
within $\Lin(H\*l_2)$, which is isomorphic to $\Cs\*p$ (see the second
paragraph of this section).

  Since $\pi\phi$ sends each $b_t$ to $b_t\*p$, we see that $\pi\phi$ is
the canonical isomorphism between $\Cs$ and $\Cs\*p$.  This shows that
$\phi$ is injective and hence an isomorphism onto $A$, concluding the
proof.
  \proofend

  \state Proposition \label{\Kamena}
  If $\B$ is a Fell bundle such that $\B\*\K$ is amenable, then $\B$,
itself, is amenable.

  \proof
  Consider the diagram
  $$
  \matrix{
  \Cs & {\buildrel \phi \over \rightarrow} & C^*(\B\*\K) \cr\cr
  E\downarrow\ \ &&\ \downarrow \tilde E\cr\cr
  B_e & \rightarrow & B_e\*\K \cr
  }
  $$
  where $\phi$ is the injective map described in the proof of
\lcite{\Kiso}, the vertical arrows represent the corresponding
conditional expectations, and the unlabeled horizontal arrow maps each
$b_e$ to $b_e\*p$.  It is easy to see that this is a commutative
diagram.  Recall that a Fell bundle is amenable if and only if the
conditional expectation on its cross-sectional algebra is faithful, as
seen in \lcite{\FaithfullCondExp}.  So, let us show that $\B$ possesses
this property.  If $x\in\Cs$ is such that $E(x^*x)=0$, then we have that
$\tilde E(\phi(x^*x))=0$, whence $\phi(x^*x)=0$ and, finally, $x=0$
because $\phi$ is injective.
  \proofend

  The following is our main result:

  \state Theorem \label{\Main}
  Let $\B$ be an orthogonal, \ss, separable Fell bundle over $\F$.  Then
$\B$ is amenable.

  \proof
  All of the properties assumed on $\B$ are clearly inherited by
$\B\*\K$.  In addition, the latter possesses a stable unit fiber algebra
and hence, by
  \lcite{\StableIsAmena}, it is amenable.  The conclusion then follows
from
  \lcite{\Kamena}.
  \proofend

  Some of the most useful facts about amenable Fell bundles
  (see e.g. \scite{\Amena}{4.10})
  require not only that the bundle be amenable, but that the \ap holds.
This leads one to ask whether the result above could be strengthened by
replacing amenability with the \ap.  We do not have a satisfactory
answer to this question but then, again, we do not know of any example
of an amenable Fell bundle which does not satisfy the \ap.

  \bigbreak
  \centerline{\tensc References}

  \begingroup
  \catcode`\@=0
  \nobreak\medskip\frenchspacing

@Article{\Claire,
  auth = {C. Anantharaman-Delaroche},
  title = {Systemes dynamiques non commutatifs et moyenabi\-li\-t\'e},
  journal = {Math. Ann.},
  year = {1987},
  volume = {279},
  pages = {297--315},
  NULL = {},
  author = {Claire Anantharaman-Delaroche},
  }

@Article{\CK,
  auth = {J. Cuntz and W. Krieger},
  title = {A Class of C*-algebras and Topological Markov Chains},
  journal = {Inventiones Math.},
  year = {1980},
  volume = {56},
  pages = {251--268},
  NULL = {},
  author = {Joachim Cuntz and W. Krieger},
  }

@Article{\BGR,
  auth = {L. G. Brown, P. Green and M. A. Rieffel},
  title = {Stable isomorphism and strong Morita equivalence of
C*-algebras},
  journal = {Pacific J. Math.},
  year = {1977},
  volume = {71},
  pages = {349--363},
  NULL = {},
  author = {L. G. Brown, P. Green and Marc A. Rieffel},
  }

@Article{\TPA,
  auth = {R. Exel},
  title = {Twisted Partial Actions, A Classification of Regular
{C*}-Algebraic Bundles},
  journal = {Proc. London Math. Soc.},
  year = {1997},
  volume = {74},
  pages = {417-443},
  NULL = {},
  akey = {Exel1997b},
  author = {Ruy Exel},
  number = {3},
  atrib = {IR},
  }

@Unpublished{\Inverse,
  auth = {R. Exel},
  title = {Partial Actions of Groups and Actions of Inverse Semigroups},
  institution = {Universidade de S\~ao Paulo},
  year = {1995},
  type = {preprint},
  note = {to appear in {\sl Proc. Amer. Math. Soc\/}},
  NULL = {},
  akey = {Exel1995b},
  author = {Ruy Exel},
  atrib = {S},
  }

@Unpublished{\Amena,
  auth = {R. Exel},
  title = {Amenability for {F}ell Bundles},
  institution = {Universidade de S\~ao Paulo},
  year = {1996},
  type = {preprint},
  note = {to appear in {\sl J. Reine Angew. Math\/}},
  NULL = {},
  akey = {Exel1996b},
  author = {Ruy Exel},
  atrib = {S},
  }

@Book{\FD,
  auth = {J. M. G. Fell and R. S. Doran},
  title = {Representations of *-algebras, locally compact groups, and
Banach *-algebraic bundles},
  publisher = {Academic Press},
  year = {1988},
  volume = {125 and 126},
  series = {Pure and Applied Mathematics},
  NULL = {},
  author = {J. M. G. Fell and R. S. Doran},
  }

@Book{\Greenleaf,
  auth = {F. P. Greenleaf},
  title = {Invariant means on topological groups},
  publisher = {van Nostrand-Reinhold},
  year = {1969},
  volume = {16},
  series = {Mathematical Studies},
  NULL = {},
  author = {F. P. Greenleaf},
  }

@TechReport{\Ng,
  auth = {C.-K. Ng},
  title = {Reduced Cross-sectional C*-algebras of C*-algebraic bundles
and Coactions},
  institution = {Oxford University},
  year = {1996},
  type = {preprint},
  note = {},
  NULL = {},
  author = {Chi-Keung Ng},
  }

@Article{\Nica,
  auth = {A. Nica},
  title = {C*-algebras generated by isometries and Wiener-Hopf
operators},
  journal = {J. Operator Theory},
  year = {1991},
  volume = {27},
  pages = {1--37},
  NULL = {},
  author = {Andu Nica},
  }

@Book{\Ped,
  auth = {G. K. Pedersen},
  title = {C*-Algebras and their automorphism groups},
  publisher = {Acad. Press},
  year = {1979},
  volume = {},
  series = {},
  NULL = {},
  author = {Gert K. Pedersen},
  }

@MastersThesis{\Nandor,
  auth = {N. Sieben},
  title = {C*-crossed products by partial actions and actions of inverse
semigroups},
  school = {Arizona State University},
  year = {1994},
  type = {Masters Thesis},
  NULL = {},
  author = {N\'andor Sieben},
  }

@Article{\Zettl,
  auth = {H. Zettl},
  title = {A characterization of ternary rings of operators},
  journal = {Adv. Math.},
  year = {1983},
  volume = {48},
  pages = {117--143},
  NULL = {},
  author = {H. Zettl},
  }

  \catcode`\@=12
  \endgroup

  \begingroup
  \parindent=0pt
  \bigskip
  \obeylines
  Departamento de Matem\'atica
  Universidade de S\~ao Paulo
  Rua do Mat\~ao, 1010
  05508-900 S\~ao Paulo -- Brazil
  exel@ime.usp.br
  \endgroup
  \bye